\documentclass[a4paper]{ifacconf}

\usepackage{graphicx}      % include this line if your document contains figures
\usepackage{natbib}        % required for bibliography
\usepackage[utf8]{inputenc}
\usepackage{amsmath}
\usepackage{amsfonts}
\usepackage{amssymb}
\usepackage{bbm}
\usepackage{mathtools}
\usepackage{float}
\usepackage{setspace}
\usepackage{graphicx}
\usepackage{enumitem}
\usepackage{tikz}
\usetikzlibrary{arrows,positioning,decorations.pathreplacing,shapes,arrows,positioning}
\usepackage{algpseudocode,algorithm,algorithmicx}

\algrenewcommand\algorithmicrequire{\textbf{Input:}}
\algrenewcommand\algorithmicensure{\textbf{Output:}}
\tikzstyle{bball} = [circle,shading=ball, ball color=black!100!white,
minimum size=0.4cm]
\tikzstyle{wball} = [circle,shading=ball, ball color=white!100!black,
minimum size=0.4cm]
\DeclareMathOperator{\rank}{rank}

\DeclareMathOperator{\adj}{adj}

\begin{document}
\begin{frontmatter}

\title{Topological Conditions for Identifiability of Dynamical Networks with Partial Node Measurements} 

%\thanks[footnoteinfo]{Sponsor and financial support acknowledgment
%goes here. Paper titles should be written in uppercase and lowercase
%letters, not all uppercase.}

\author[First]{Henk J. van Waarde} 
\author[Second]{Pietro Tesi} 
\author[First]{M. Kanat Camlibel}

\address[First]{Bernoulli Institute for Mathematics, Computer Science and Artificial Intelligence, University of Groningen, The Netherlands. (e-mail: h.j.van.waarde@rug.nl and m.k.camlibel@rug.nl).}
\address[Second]{Engineering and Technology Institute Groningen, University of Groningen, The Netherlands. Pietro Tesi is also with the Department of Information Engineering, University of Florence, Italy (e-mail: p.tesi@rug.nl, pietro.tesi@unifi.it).}

\begin{abstract}                % Abstract of not more than 250 words.
This paper deals with dynamical networks for which the relations between node signals are described by proper transfer functions and external signals can influence each of the node signals. In particular, we are interested in graph-theoretic conditions for \emph{identifiability} of such dynamical networks, where we assume that only a subset of nodes is measured but the underlying graph structure of the network is known. This problem has recently been investigated in the case of \emph{generic} identifiability. In this paper, we investigate a stronger notion of identifiability \emph{for all} network matrices associated with a given graph. For this, we introduce a new graph-theoretic concept called \emph{constrained vertex-disjoint paths}. As our main result, we state conditions for identifiability based on these constrained vertex-disjoint paths.
\end{abstract}


\end{frontmatter}
%===============================================================================
\section{Introduction}
Networks of dynamical systems appear in a variety of domains, including power systems, robotic networks and water distribution networks. In this paper, we consider dynamical networks where the relations between node signals are modelled by proper transfer functions and external signals can influence each of the node signals. Such network models have received much attention (see, e.g., [\cite{vandenHof2013}], [\cite{Dankers2014}], [\cite{Hendrickx2018}], [\cite{Weerts2018}]). The interconnection structure of a dynamical network can be represented by a directed graph, where vertices (or nodes) represent scalar signals, and edges correspond to transfer functions connecting different node signals. We will assume that the underlying graph (i.e., the topology) of the dynamical network is \emph{known}. We remark that identification of the network topology has also been studied (see, e.g., [\cite{Goncalves2008}], [\cite{Shahrampour2015}], [\cite{Hayden2017}], [\cite{vanWaarde2017}]).

We are interested in conditions for \emph{identifiability} of dynamical networks. Loosely speaking, identifiability comprises the ability to distinguish between certain (network) models on the basis of measured data. In this work, we assume that each node of the network is externally excited by a known signal. However, the node signals of only a subset of nodes is measured. Within this setup, we are interested in two identifiability problems. Firstly, we want to find conditions under which the transfer functions from a given node to its out-neighbours can be identified. Secondly, we wonder under which conditions \emph{all} transfer functions in the network can be identified. In particular, our aim is to find \emph{graph-theoretic} conditions for the above problems. Conditions based on the topology of the network are desirable, since they give insight in the types of network structures that allow identification, and in addition may aid in the \emph{selection} of measured nodes.  

In previous work [\cite{Hendrickx2018}], graph-theoretic conditions have been established for \emph{generic} identifiability. That is, conditions were given under which a subset of transfer functions in the network can be identified for \emph{``almost all"} network matrices associated with a given graph. In contrast to [\cite{Hendrickx2018}], we are interested in graph-theoretic conditions for a stronger notion, namely identifiability for \emph{all} network matrices associated with the graph. The difference between generic identifiability and identifiability for all network matrices might seem subtle at first, however, similar differences in the \emph{controllability} literature have led to very different graph-theoretic characterizations (cf. [\cite{Liu2011}], [\cite{Chapman2013}]).

In this paper, we introduce a new graph-theoretic concept called \emph{constrained vertex-disjoint paths}. Such paths are a generalization of \emph{constrained matchings} in bipartite graphs [\cite{Hershkowitz1993}]. As our main result, we state conditions for identifiability in terms of constrained vertex-disjoint paths. 

This paper is organized as follows. In Section \ref{sectionpreliminaries} we introduce notation and preliminaries. Subsequently, in Section \ref{sectionproblem} we state the problem. Next, Section \ref{sectionresults} contains the main results, and in Section \ref{sectionconclusions} we give our conclusions.

\section{Preliminaries}
\label{sectionpreliminaries}
We denote the set of natural numbers by $\mathbb{N}$, real numbers by $\mathbb{R}$, complex numbers by $\mathbb{C}$, and real $m \times n$ matrices by $\mathbb{R}^{m \times n}$. The $n \times n$ identity matrix is denoted by $I_n$. When its dimension is clear from the context, we simply write $I$. Consider a rational function $f(z) = \frac{p(z)}{q(z)}$, where $p(z)$ and $q(z)$ are polynomials with real coefficients. Then, $f$ is called \emph{proper} if the degree of $p(z)$ is less than or equal to the degree of $q(z)$. We say $f$ is \emph{strictly proper} if the degree of $p(z)$ is less than the degree of $q(z)$. An $m \times n$ matrix $A(z)$ is called \emph{rational} if its entries are rational functions in the indeterminate $z$. In addition, $A(z)$ is \emph{proper} if its entries are proper rational functions in $z$. We omit the argument $(z)$ whenever the dependency of $A$ on $z$ is clear from the context. The \emph{normal rank} of $A(z)$ is defined as $\max_{\lambda \in \mathbb{C}} \rank A(\lambda)$ and denoted by $\rank A(z)$, with slight abuse of notation. We denote the $(i,j)$-th entry of $A$ by $A_{ij}$. Moreover, the $j$-th column of $A$ is given by $A_{\bullet j}$. More generally, let $\mathcal{M} \subseteq \{1,2,\dots,m\}$ and $\mathcal{N} \subseteq \{1,2,\dots,n\}$. Then, $A_{\mathcal{M},\mathcal{N}}$ denotes the submatrix of $A$ containing the rows of $A$ indexed by $\mathcal{M}$ and the columns of $A$ indexed by $\mathcal{N}$. Next, consider the case that $A$ is square, i.e., $m = n$. The \emph{determinant} of $A$ is denoted by $\det A$, while the \emph{adjugate} of $A$ is denoted by $\adj A$. A \emph{principal submatrix} of $A$ is a submatrix $A_{\mathcal{M},\mathcal{M}}$, where $\mathcal{M} \subseteq \{1,2,\dots,m\}$. The determinant of $A_{\mathcal{M},\mathcal{M}}$ is called a \emph{principal minor} of $A$. 

\subsection{Graph theory}
\label{sectiongraphtheory}
Let $\mathcal{G} = (\mathcal{V},\mathcal{E})$ be a directed graph, with vertex set $\mathcal{V} = \{1,2,\dots,n\}$ and edge set $\mathcal{E} \subseteq \mathcal{V} \times \mathcal{V}$. Unless explicitly mentioned, the graphs considered in this paper are \emph{simple}, i.e., without self-loops (edges of the form $(i,i)\in\mathcal{E}$) and with at most one edge from one node to another. A node $j \in \mathcal{V}$ is said to be an \emph{out-neighbour} of $i \in \mathcal{V}$ if $(i,j) \in \mathcal{E}$. The set of out-neighbours of node $i \in \mathcal{V}$ is denoted by $\mathcal{N}_i$. For any subset $\mathcal{S} = \{v_1,v_2,\dots,v_s\} \subseteq \mathcal{V}$ we define the $s \times n$ matrix $P(\mathcal{V};\mathcal{S})$ as $P_{ij} := 1$ if $j = v_i$, and $P_{ij} := 0$ otherwise. The complement of $\mathcal{S}$ in $\mathcal{V}$ is defined as $\mathcal{S}^c := \mathcal{V} \setminus \mathcal{S}$. Moreover, the cardinality of $\mathcal{S}$ is denoted by $|\mathcal{S}|$. A \emph{path} $\mathcal{P}$ is a set of edges in $\mathcal{G}$ of the form $\mathcal{P} = \{(v_i,v_{i+1}) \mid i = 1,2,\dots,k \} \subseteq \mathcal{E}$, where the vertices $v_1,v_2,\dots,v_{k+1}$ are \emph{distinct}. The vertex $v_1$ is called a \emph{starting node} of $\mathcal{P}$, while $v_{k+1}$ is the \emph{end node}. The cardinality of $\mathcal{P}$ is called the \emph{length} of the path.  A collection of paths $\mathcal{P}_1,\mathcal{P}_2,\dots,\mathcal{P}_l$ is called \emph{vertex-disjoint} if the paths have no vertex in common, that is, if for all distinct $i,j \in \{1,2,\dots,l\}$, we have that
\begin{equation*}
(u_i,w_i) \in \mathcal{P}_i, (u_j,w_j) \in \mathcal{P}_j \implies u_i,w_i,u_j,w_j \text{ are distinct}.
\end{equation*}
Consider two disjoint subsets $\mathcal{V}_1, \mathcal{V}_2 \subseteq \mathcal{V}$. We say there exist $m$ vertex-disjoint paths \emph{from} $\mathcal{V}_1$ \emph{to} $\mathcal{V}_2$ if there exist $m$ vertex-disjoint paths in $\mathcal{G}$ with starting nodes in $\mathcal{V}_1$ and end nodes in $\mathcal{V}_2$. In the case that $\mathcal{V}_1 \cap \mathcal{V}_2 \neq \emptyset$, we say that there exist $m$ vertex-disjoint paths \emph{from} $\mathcal{V}_1$ \emph{to} $\mathcal{V}_2$ if there are $\max\{0,m-|\mathcal{V}_1 \cap \mathcal{V}_2|\}$ vertex-disjoint paths from $\mathcal{V}_1 \setminus \mathcal{V}_2$ to $\mathcal{V}_2 \setminus \mathcal{V}_1$. A \emph{cycle} $\mathcal{K}$ is a set of edges in $\mathcal{G}$ of the form $\mathcal{K} = \{(v_i,v_{i+1}) \mid i = 1,2,\dots,k \} \subseteq \mathcal{E}$, where $v_1,v_2,\dots,v_{k}$ are distinct, and $v_1 = v_{k+1}$. The cardinality of $\mathcal{K}$ is called the \emph{length} of the cycle. A collection of cycles is called \emph{vertex-disjoint} if the cycles have no vertex in common. Moreover, a \emph{spanning cycle family} in $\mathcal{G}$ is a collection of vertex-disjoint cycles such that each vertex in $\mathcal{V}$ is contained in exactly one cycle. Next, consider a \emph{weighted} directed graph $\mathcal{G} = (\mathcal{V},\mathcal{E})$, that is, a directed graph where each edge $(i,j) \in \mathcal{E}$ has an associated rational function $f_{ji}(z)$ called the \emph{weight} of the edge $(i,j)$. The weight of a path $\mathcal{P}$ in $\mathcal{G}$ is defined as the product of the weights of all edges in $\mathcal{P}$. Moreover, the weight of a set of vertex-disjoint paths is defined as the product of the weights of all paths in the set. Similarly, the weight of a cycle $\mathcal{K}$ is defined as the product of the weights of all edges in $\mathcal{K}$. Finally, the weight of a set of vertex-disjoint cycles is defined as the product of weights of all cycles in the set.

\section{Problem statement and motivation}
\label{sectionproblem}
Let $\mathcal{G} = (\mathcal{V},\mathcal{E})$ be a simple directed graph, with vertex set $\mathcal{V} = \{1,2,\dots,n\}$ and edge set $\mathcal{E} \subseteq \mathcal{V} \times \mathcal{V}$. Following the setup of [\cite{Hendrickx2018}] (see also [\cite{vandenHof2013}], [\cite{Weerts2018}]), we associate the following dynamical system with the graph $\mathcal{G}$:
\begin{equation}
\label{system}
\begin{aligned}
w(t) &= G(q) w(t) + r(t) + v(t) \\
y(t) &= C w(t).
\end{aligned} 
\end{equation}
Here $w, r$, and $v$ are $n$-dimensional vectors of node signals, known external signals, and unknown disturbances, respectively. The output vector $y$ is $p$-dimensional, and consists of the node signals of a subset $\mathcal{C} \subseteq \mathcal{V}$ of nodes, with $|\mathcal{C}| = p$. Consequently, the matrix $C$ is defined as $C := P(\mathcal{V},\mathcal{C})$. Moreover, $q^{-1}$ denotes the unit delay operator, i.e., $q^{-1} w_i(t) = w_i(t-1)$. Finally, $G(z)$ is an $n \times n$ rational matrix, called the \emph{network matrix}, satisfying the following properties [\cite{vandenHof2013}]:
\begin{enumerate}[label=\textbf{P\arabic*.}]
	\item For all $i,j \in \mathcal{V}$, the entry $G_{ji}(z)$ is a proper rational (transfer) function.
	\item The function $G_{ji}(z)$ is nonzero if and only if $(i,j) \in \mathcal{E}$. A matrix $G(z)$ that satisfies this property is said to be \emph{consistent} with the graph $\mathcal{G}$.
	\item Every principal minor of $\lim_{z \to \infty} (I-G(z))$ is nonzero. This implies that the network model \eqref{system} is \emph{well-posed} in the sense of Definition 2.11 of [\cite{Dankers2014}].
\end{enumerate}
Property P3 is required for the technical analysis in this paper, but only imposes weak restrictions on the matrix $G$ (see [\cite{vandenHof2013}]). 
\begin{rem}
We focus on the network model \eqref{system} that was originally introduced in [\cite{vandenHof2013}]. Note that state-space network models have also received much attention (see, e.g., [\cite{Goncalves2008}], [\cite{Hayden2017}]). A state-space model (with scalar node dynamics) can be obtained from \eqref{system} by choosing the nonzero entries of $G$ as first-order strictly proper functions. However, the model \eqref{system} is more general in the sense that higher-order transfer functions are also allowed. 
\end{rem}
A network matrix $G$ satisfying Properties P1, P2, and P3 is called \emph{admissible}. In what follows, we use the shorthand notation $T(z) := (I - G(z))^{-1}$, where $G$ is assumed to be admissible. Note that \eqref{system} implies that $y(t) = CT(q) r(t) + CT(q) v(t)$, which shows that the transfer matrix from $r$ to $y$ is given by $CT(z)$. In this paper, we are interested in the question of which transfer functions in $G(z)$ can be uniquely identified from input/output data, that is, from the external signals $r(t)$ and output signals $y(t)$. To this end, we assume that the graph $\mathcal{G} = (\mathcal{V},\mathcal{E})$ is \emph{known}. Moreover, we assume that the excitation signal $r(t)$ is sufficiently rich such that, under suitable assumptions on the disturbance $v(t)$, the transfer matrix $CT(z)$ can be identified from $\{r(t),y(t)\}$-data [\cite{Ljung1999}]. Note that we are not per se interested in identifying the matrix $CT(z)$, but we want to find (a part of) the network matrix $G(z)$. Therefore, the question is which transfer functions in $G(z)$ can be uniquely reconstructed from the transfer matrix $CT(z)$. In recent work [\cite{Hendrickx2018}], this question has been considered for \emph{generic identifiability}. Graph-theoretic conditions were given under which a set of transfer functions can be uniquely identified from $CT(z)$ \emph{for almost all} network matrices $G(z)$ consistent with the graph $\mathcal{G}$. For a formal definition of generic identifiability we refer to Definition 1 of [\cite{Hendrickx2018}]. Here, we will informally illustrate the approach of [\cite{Hendrickx2018}]. 
\begin{exmp}
	\label{example1}
	Consider the graph $\mathcal{G} = (\mathcal{V},\mathcal{E})$ depicted in Figure \ref{fig:graph1}. We assume that the node signals of nodes $4$ and $5$ can be measured, that is, $\mathcal{C} = \{4,5\}$. Suppose that we want to identify the transfer functions from node $1$ to its out-neighbours, i.e., the transfer functions $G_{21}(z)$ and $G_{31}(z)$. According to Corollary 5.1 of [\cite{Hendrickx2018}], this is possible if and only if there exist two vertex-disjoint paths from $\mathcal{N}_1$ to $\mathcal{C}$. Note that this is the case in this example since the edges $(2,4)$ and $(3,5)$ are two vertex-disjoint paths. 
	\begin{figure}[h!]
		\centering
		\scalebox{0.8}{
			\begin{tikzpicture}[scale=1]
			\node[draw, style=wball, label={180:$1$}] (1) at (0,0) {};
			\node[draw, style=wball, label={90:$2$}] (2) at (2,1) {};
			\node[draw, style=wball, label={-90:$3$}] (3) at (2,-1) {};
			\node[draw, style=wball, label={90:$4$}] (4) at (4,1) {};
			\node[draw, style=wball, label={-90:$5$}] (5) at (4,-1) {};
			\draw[-latex] (1) -- (2);
			\draw[-latex] (1) -- (3);
			\draw[-latex] (2) -- (4);
			\draw[-latex] (2) -- (5);
			\draw[-latex] (3) -- (4);
			\draw[-latex] (3) -- (5);
			\end{tikzpicture}
		}
		\caption{Graph used in Example \ref{example1}.}
		\label{fig:graph1}
	\end{figure}
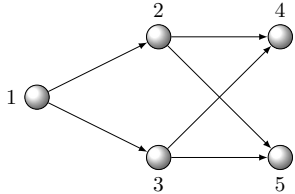
	To see why we can generically identify the transfer functions $G_{21}$ and $G_{31}$, we compute the transfer matrix $CT$ as follows:
	\begin{equation*}
	CT = \begin{pmatrix}
	G_{42} G_{21} + G_{43} G_{31} & G_{42} & G_{43} & 1 & 0 \\
	G_{52} G_{21} + G_{53} G_{31} & G_{52} & G_{53} & 0 & 1
	\end{pmatrix}.
	\end{equation*}
	Clearly, we can uniquely identify the transfer functions $G_{42}, G_{43}, G_{52}$, and $G_{53}$ from $CT$. Moreover, the transfer matrices $G_{21}$ and $G_{31}$ satisfy 
	\begin{equation}
	\label{equationex1}
	\begin{pmatrix}
	G_{42} & G_{43} \\
	G_{52} & G_{53}
	\end{pmatrix}
	\begin{pmatrix}
	G_{21} \\ G_{31}
	\end{pmatrix}
	= \begin{pmatrix}
	T_{41} \\ T_{51}
	\end{pmatrix}.
	\end{equation} 
	Equation \eqref{equationex1} has a unique solution in the unknowns $G_{21}$ and $G_{31}$ if $G_{42} G_{53} - G_{43} G_{52} \neq 0$, which means that we can identify $G_{21}$ and $G_{31}$ for almost all $G$ consistent with the graph $\mathcal{G}$ (specifically, for all $G$ except those for which $G_{42} G_{53} - G_{43} G_{52} = 0$). $\hfill{} \blacksquare$
\end{exmp}
As mentioned before, the approach based on vertex-disjoint paths [\cite{Hendrickx2018}] gives necessary and sufficient conditions for \emph{generic} identifiability. This implies that for some network matrices $G$, it might be impossible to identify the transfer functions, even though the path-based conditions are satisfied. For instance, in Example \ref{example1} we cannot identify the transfer functions $G_{21}$ and $G_{31}$ if the network matrix $G$ is such that $G_{42} = G_{43} = G_{52} = G_{53}$. Nonetheless, a scenario in which some of the transfer matrices in the network are equal may occur in practice.

Instead of generic identifiability, in this paper we are interested in graph-theoretic conditions that guarantee identifiability \emph{for all} network matrices consistent with a given graph. Such a problem might seem like a simple extension of the work on generic identifiability [\cite{Hendrickx2018}]. However, to analyze \emph{strong structural} network properties (for all network matrices), we typically need different graph-theoretic tools than the ones used in the analysis of \emph{generic} network properties. For instance, in the literature on \emph{controllability} of dynamical networks, generic controllability (often called (weak) structural controllability) is related to maximal matchings [\cite{Liu2011}], while strong structural controllability is related to constrained matchings [\cite{Chapman2013}]. To make the problem of this paper more precise, we state a few definitions. Firstly, we are interested in conditions under which all transfer functions from a node $i$ to its out-neighbours $\mathcal{N}_i$ are identifiable (for any admissible network matrix $G$, i.e., any $G$ that satisfies properties P1, P2, and P3). If this is the case, we say $(i,\mathcal{N}_i)$ is identifiable. More precisely, we have the following definition.

\begin{defn}
	Consider a directed graph $\mathcal{G} = (\mathcal{V},\mathcal{E})$. Let $i \in \mathcal{V}$, $\mathcal{C} \subseteq \mathcal{V}$, and define $C := P(\mathcal{V},\mathcal{C})$. We say $(i,\mathcal{N}_i)$ is \emph{identifiable} from $\mathcal{C}$ if the implication
	\begin{equation*}
	C(I - G(z))^{-1} = C(I - \bar{G}(z))^{-1} \implies G_{\bullet i}(z) = \bar{G}_{\bullet i}(z)
	\end{equation*}
	holds for all admissible network matrices $G(z)$ and $\bar{G}(z)$.
\end{defn}

In addition to identifiability of $(i,\mathcal{N}_i)$, we are interested in conditions under which the entire network matrix $G$ can be identified from the transfer matrix $CT$. If this is the case, we say the graph $\mathcal{G}$ is identifiable.

\begin{defn}
	Consider a directed graph $\mathcal{G} = (\mathcal{V},\mathcal{E})$. Let $\mathcal{C} \subseteq \mathcal{V}$ and define $C := P(\mathcal{V},\mathcal{C})$. We say $\mathcal{G}$ is \emph{identifiable} from $\mathcal{C}$ if the implication
	\begin{equation*}
	C(I - G(z))^{-1} = C(I - \bar{G}(z))^{-1} \implies G(z) = \bar{G}(z)
	\end{equation*}
	holds for all admissible network matrices $G(z)$ and $\bar{G}(z)$.
\end{defn}

The main goal of this paper is to find graph-theoretic conditions under which $(i,\mathcal{N}_i)$ is identifiable. Furthermore, based on such conditions, we want to establish graph-theoretic conditions under which $\mathcal{G}$ is identifiable.

\section{Main results}
\label{sectionresults}

In this section we will present our main results. In Section \ref{subsectionrank} we give conditions for necessary and sufficient rank conditions for identifiability. Subsequently, in Section \ref{subsectiongraph} we use these rank conditions to derive graph-theoretic tests for identifiability.  

\subsection{Rank conditions for identifiability}
\label{subsectionrank}
First, we give a condition for identifiability of $(i,\mathcal{N}_i)$ in terms of the normal rank of $T_{\mathcal{C}, \mathcal{N}_i}(z)$ in Lemma \ref{lemma1}. Note that this condition is similar to the one in Lemma 5.1 of [\cite{Hendrickx2018}], however, since we have additional assumptions (P2 and P3) on the network matrix, the result of [\cite{Hendrickx2018}] is not directly applicable to our setup. Therefore, we provide a proof of Lemma \ref{lemma1}.
\begin{lem}
	\label{lemma1}
	Consider a directed graph $\mathcal{G} = (\mathcal{V},\mathcal{E})$. Let $i \in \mathcal{V}$ and $\mathcal{C} \subseteq \mathcal{V}$. Then, $(i,\mathcal{N}_i)$ is identifiable from $\mathcal{C}$ if and only if for any admissible network matrix $G(z)$ we have $\rank T_{\mathcal{C}, \mathcal{N}_i}(z) = | \mathcal{N}_i |$, where $T(z) := (I-G(z))^{-1}$.
\end{lem}

\begin{pf}
	For the `if' part, suppose that $\rank T_{\mathcal{C}, \mathcal{N}_i}(z) = | \mathcal{N}_i |$ for any admissible network matrix $G(z)$. Let $\bar{G}(z)$ be an admissible network matrix satisfying $C(I-\bar{G}(z))^{-1} = C(I-G(z))^{-1}$. We define $D(z) := G(z) - \bar{G}(z)$, and note that the following four statements are equivalent:
	\begin{equation}
	\label{derivationrank}
	\begin{aligned}
	C &= C(I-G(z))^{-1} (I-\bar{G}(z)) \\
	C &= C(I-G(z))^{-1} (I-G(z) + D(z)) \\
	C &= C + C(I-G(z))^{-1} D(z) \\
	0 &= CT(z) D(z).
	\end{aligned}
	\end{equation}
	In particular, we obtain $CT(z) D_{\bullet i}(z) = 0$. Since both $G$ and $\bar{G}$ are consistent with the graph, we have that $D_{ji}(z) = 0$ if $j \not\in \mathcal{N}_i$. Consequently, $T_{\mathcal{C}, \mathcal{N}_i}(z) D_{\mathcal{N}_i, i}(z) = 0$. By hypothesis, $\rank T_{\mathcal{C}, \mathcal{N}_i}(z) = | \mathcal{N}_i |$, and therefore $D_{\mathcal{N}_i, i}(z) = 0$. Hence, also $D_{\bullet i}(z) = 0$. We conclude that $G_{\bullet i}(z) = \bar{G}_{\bullet i}(z)$, which shows that $(i,\mathcal{N}_i)$ is identifiable. 
	
	For the `only if' part, suppose that $\rank T_{\mathcal{C}, \mathcal{N}_i}(z) < | \mathcal{N}_i |$ for some admissible network matrix $G$. We want to prove that $(i,\mathcal{N}_i)$ is not identifiable, that is, we want to prove the existence of an admissible network matrix $\bar{G}$ such that $C(I - G(z))^{-1} = C(I - \bar{G}(z))^{-1}$, but $G_{\bullet i}(z) \neq \bar{G}_{\bullet i}(z)$. Note that by our hypothesis, there exists a nonzero rational vector $\hat{w}(z)$ such that $T_{\mathcal{C}, \mathcal{N}_i}(z) \hat{w}(z) = 0$. Obviously, this means that $T_{\mathcal{C}, \mathcal{N}_i}(z) \hat{v}(z) = 0$, where $\hat{v}(z) := \alpha z^{-k} \hat{w}(z)$ for any $k \in \mathbb{N}$ and $\alpha \in \mathbb{R}$. We choose $k \in \mathbb{N}$ such that the nonzero entries of $\hat{v}(z)$ are strictly proper. Moreover, we choose $\alpha \in \mathbb{R} \setminus \{ 0 \}$ such that all entries of $\hat{G}_{\bullet i}(z) - \hat{v}(z)$ are nonzero, where $\hat{G}_{\bullet i}(z)$ is the vector obtained from $G_{\bullet i}(z)$ by removing the entries corresponding to $\mathcal{V} \setminus \mathcal{N}_i$. Note that this is always possible, since $G$ is consistent with the graph, and therefore all entries of $\hat{G}_{\bullet i}(z)$ are nonzero. Let $v(z)$ denote the $n$-dimensional rational vector with the following two properties. Firstly, $v(z)$ has zero entries in positions corresponding to nodes in $\mathcal{V} \setminus \mathcal{N}_i$. Secondly, the vector obtained from $v(z)$ by removing all entries corresponding to $\mathcal{V} \setminus \mathcal{N}_i$ equals $\hat{v}(z)$. In addition, let the versor $u \in \mathbb{R}^n$ be such that $u_i = 1$ and $u_j = 0$ for all $j \in \mathcal{V} \setminus \{i\}$. We define the matrix $\bar{G}(z) := G(z) - v(z) u^\top$. Moreover, define $D(z) := v(z) u^\top$, and note that \eqref{derivationrank} yields $C(I - G(z))^{-1} = C(I - \bar{G}(z))^{-1}$. Furthermore, since $v(z)$ is nonzero, we immediately obtain $G_{\bullet i}(z) \neq \bar{G}_{\bullet i}(z)$. It remains to be shown that the matrix $\bar{G}$ is admissible, i.e., that $\bar{G}$ satisfies properties P1, P2, and P3. 
	
	Firstly, note that $\bar{G}(z) = G(z) - v(z) u^\top$ is the difference of a proper transfer matrix $G(z)$ and a matrix $v(z) u^\top$ with entries that are either zero or strictly proper transfer functions. Consequently, each entry of $\bar{G}(z)$ is a proper rational function, and $\bar{G}$ satisfies P1. Secondly, to prove that $\bar{G}$ satisfies P2, note that for all $k \in \mathcal{V} \setminus \{i\}$ and all $j \in \mathcal{V}$ we have $\bar{G}_{jk}(z) = G_{jk}(z)$, and consequently, $\bar{G}_{jk}(z) \neq 0$ if and only if $(k,j) \in \mathcal{E}$. Moreover, by construction of $v(z)$, we have that $\bar{G}_{ji}(z) \neq 0$ if and only if $j \in \mathcal{N}_i$. We conclude that $\bar{G}$ is consistent with the graph. Finally, to prove that $\bar{G}$ satisfies property P3, note that 
	\begin{equation}
	\label{sumlimits}
	\lim_{z \to \infty} (I - \bar{G}(z)) = \lim_{z \to \infty} (I-G(z)) + \lim_{z \to \infty} v(z) u^\top,
	\end{equation} 
	since both limits on the right-hand side of \eqref{sumlimits} exist. In fact, since an entry of $v(z)$ is either zero or strictly proper, we have $\lim_{z \to \infty} v(z) u^\top = 0$, and consequently 
	\begin{equation*}
	\lim_{z \to \infty} (I - \bar{G}(z)) = \lim_{z \to \infty} (I-G(z)).
	\end{equation*} 
	This shows that $\bar{G}$ satisfies property P3 (as $G$ satisfies P3). To conclude, we have shown the existence of an admissible $\bar{G}$ such that $C(I - G(z))^{-1} = C(I - \bar{G}(z))^{-1}$, but $G_{\bullet i}(z) \neq \bar{G}_{\bullet i}(z)$. That is, $(i,\mathcal{N}_i)$ is not identifiable. This proves the lemma. $\hfill{} \square$
\end{pf}

As an immediate consequence of Lemma \ref{lemma1}, we find conditions for the identifiability of $\mathcal{G}$. 
\begin{cor}
	Consider a directed graph $\mathcal{G} = (\mathcal{V},\mathcal{E})$ and let $\mathcal{C} \subseteq \mathcal{V}$. Then, $\mathcal{G}$ is identifiable from $\mathcal{C}$ if and only if for all $i \in \mathcal{V}$ and all admissible $G(z)$ we have $\rank T_{\mathcal{C}, \mathcal{N}_i}(z) = | \mathcal{N}_i |$, where $T(z) := (I-G(z))^{-1}$.
\end{cor}

In the case that $T_{\mathcal{C}, \mathcal{N}_i}(z)$ is square, we can relate the rank condition $\rank T_{\mathcal{C}, \mathcal{N}_i}(z) = |\mathcal{N}_i|$ to a rank condition of a certain submatrix of $I - G(z)$. Specifically, the following lemma states the equivalence of $\det T_{\mathcal{C}, \mathcal{N}_i}(z) \neq 0$ and $\det\left( (I-G(z))_{\mathcal{N}_i^c, \mathcal{C}_{\phantom{i}}^c} \right) \neq 0$, where the nonzero condition should be understood as nonzero \emph{as a rational function}.

\begin{lem}
	\label{lemma2}
	Consider a directed graph $\mathcal{G} = (\mathcal{V},\mathcal{E})$ and let $i \in \mathcal{V}$. Suppose that $\mathcal{C} \subseteq \mathcal{V}$ satisfies $| \mathcal{C} | = | \mathcal{N}_i |$. Let $G(z)$ be an admissible network matrix and define $T(z) := (I-G(z))^{-1}$. Then, $\det T_{\mathcal{C}, \mathcal{N}_i}(z) \neq 0$ if and only if $\det\left( (I-G(z))_{\mathcal{N}_i^c, \mathcal{C}_{\phantom{i}}^c} \right) \neq 0$.
\end{lem}

\begin{pf}
	We define $A(z) := \adj(I - G(z))$. Note that
	\begin{equation*}
	\begin{aligned}
	\det T_{\mathcal{C}, \mathcal{N}_i}(z) \neq 0 &\iff \det A_{\mathcal{C}, \mathcal{N}_i}(z) \neq 0 \\
	&\iff \det \left( (A^\top(z))_{\mathcal{N}_i, \mathcal{C}} \right) \neq 0.
	\end{aligned}
	\end{equation*}
	Next, we apply Jacobi's identity for the determinant of a submatrix of the adjugate (cf. Theorem 2.5.2 of [\cite{Prasolov1994}]), which shows that $\det \left( (A^\top(z))_{\mathcal{N}_i, \mathcal{C}} \right)$ is equal to
	\begin{equation*}
	\begin{aligned}
	 \pm \det\left( (I-G(z))_{\mathcal{N}_i^c, \mathcal{C}_{\phantom{i}}^c} \right) \det(I - G(z))^{|\mathcal{C}|-1}.
	\end{aligned}
	\end{equation*}
	Since $\det(I - G(z)) \neq 0$, this proves the lemma. $\hfill{} \square$
\end{pf}

\subsection{Graph-theoretic conditions for identifiability}
\label{subsectiongraph}
In this section, we provide graph-theoretic conditions for identifiability. To give some intuition for the approach, we start with a simple example.
\begin{exmp}
\label{example3}
Consider the graph in Figure \ref{fig:graph1}. We saw in Example \ref{example1} that $(1,\mathcal{N}_1)$ is not identifiable from $\mathcal{C} = \{4,5\}$ since \eqref{equationex1} has multiple solutions in $G_{21}$ and $G_{31}$ in the case that $G_{42} G_{53} - G_{43} G_{52} = 0$. Suppose that we consider a slightly different graph, namely the one in Figure \ref{fig:graph3}.
\begin{figure}[h!]
	\centering
	\scalebox{0.8}{
		\begin{tikzpicture}[scale=1]
		\node[draw, style=wball, label={180:$1$}] (1) at (0,0) {};
		\node[draw, style=wball, label={90:$2$}] (2) at (2,1) {};
		\node[draw, style=wball, label={-90:$3$}] (3) at (2,-1) {};
		\node[draw, style=wball, label={90:$4$}] (4) at (4,1) {};
		\node[draw, style=wball, label={-90:$5$}] (5) at (4,-1) {};
		\draw[-latex] (1) -- (2);
		\draw[-latex] (1) -- (3);
		\draw[-latex] (2) -- (4);
		\draw[-latex] (3) -- (4);
		\draw[-latex] (3) -- (5);
		\end{tikzpicture}
	}
	\caption{Graph used in Example \ref{example3}.}
	\label{fig:graph3}
\end{figure}
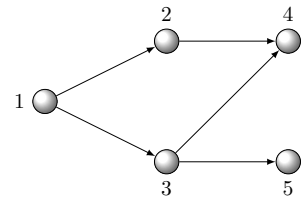
We are still interested in identifiability of $(1,\mathcal{N}_1)$ from $\mathcal{C} = \{4,5\}$. A simple calculation shows that in this case
\begin{equation}
\label{equationex3}
\begin{pmatrix}
G_{42} & G_{43} \\
0 & G_{53}
\end{pmatrix}
\begin{pmatrix}
G_{21} \\ G_{31}
\end{pmatrix}
= \begin{pmatrix}
T_{41} \\ T_{51}
\end{pmatrix},
\end{equation}
where $G_{42}, G_{43}, G_{53}$ are transfer functions that can be identified from $CT$. Note that the matrix on the left-hand side of \eqref{equationex3} has full rank for all nonzero $G_{42}, G_{43},$ and $G_{53}$. Consequently, we can identify the transfer functions $G_{21}$ and $G_{31}$ for all admissible network matrices $G$. That is, $(1,\mathcal{N}_1)$ is identifiable from $\mathcal{C}$. We observe the following difference between the graphs in Figures \ref{fig:graph1} and \ref{fig:graph3}: in Figure \ref{fig:graph1}, there are two different sets of two vertex-disjoint paths from $\{2,3\}$ to $\{4,5\}$, while in the graph of Figure \ref{fig:graph3} there exists exactly one set of two vertex-disjoint paths between these vertices. Therefore, it seems that identifiability of $(i,\mathcal{N}_i)$ does not only depend on the \emph{existence} of $|\mathcal{N}_i|$ vertex-disjoint paths from $\mathcal{N}_i$ to $\mathcal{C}$ (as is the case for generic identifiability [\cite{Hendrickx2018}]), but also depends on the \emph{number} of such sets of vertex-disjoint paths. $\hfill{} \blacksquare$
\end{exmp}

To make the idea of Example \ref{example3} more precise, we need the following definition of \emph{constrained vertex-disjoint paths}. 

\begin{defn}
	Let $\mathcal{G} = (\mathcal{V},\mathcal{E})$ be a directed graph. Consider a set of $m$ vertex-disjoint paths in $\mathcal{G}$ with starting nodes $\bar{\mathcal{V}}_1 \subseteq \mathcal{V}$ and end nodes $\bar{\mathcal{V}}_2 \subseteq \mathcal{V}$. We say that the set of vertex-disjoint paths is \emph{constrained} if it is the \emph{only} set of $m$ vertex-disjoint paths from $\bar{\mathcal{V}}_1$ to $\bar{\mathcal{V}}_2$.
\end{defn} 

Next, let $\mathcal{V}_1, \mathcal{V}_2 \subseteq \mathcal{V}$ be disjoint subsets. We say that there exists a constrained set of $m$ vertex-disjoint paths \emph{from} $\mathcal{V}_1$ \emph{to} $\mathcal{V}_2$ if there exists a constrained set of $m$ vertex-disjoint paths in $\mathcal{G}$ with starting nodes $\bar{\mathcal{V}}_1 \subseteq \mathcal{V}_1$ and end nodes $\bar{\mathcal{V}}_2 \subseteq \mathcal{V}_2$. In the case that $\mathcal{V}_1 \cap \mathcal{V}_2 \neq \emptyset$, we say that there is a constrained set of $m$ vertex-disjoint paths from $\mathcal{V}_1$ to $\mathcal{V}_2$ if there exists a constrained set of $\max\{0,m - |\mathcal{V}_1 \cap \mathcal{V}_2|\}$ vertex-disjoint paths from $\mathcal{V}_1 \setminus \mathcal{V}_2$ to $\mathcal{V}_2 \setminus \mathcal{V}_1$. Roughly speaking, this means that we count paths of ``length zero" from every node in $\mathcal{V}_1 \cap \mathcal{V}_2$ to itself.

\begin{rem}
	Note that for a set of $m$ vertex-disjoint paths from $\mathcal{V}_1$ to $\mathcal{V}_2$ to be constrained, we do not require there to be a unique set of $m$ vertex-disjoint paths from $\mathcal{V}_1$ to $\mathcal{V}_2$. In fact, we only require there to be a unique set of vertex-disjoint paths between the \emph{starting nodes} $\bar{\mathcal{V}}_1$ of the paths and the \emph{end nodes} $\bar{\mathcal{V}}_2$. 
\end{rem}

\begin{rem}
	The notion of constrained vertex-disjoint paths is strongly related to the notion of \emph{constrained matchings} in bipartite graphs [\cite{Hershkowitz1993}]. In fact, a constrained matching can be seen as a special case of a constrained set of vertex-disjoint paths, where all paths are of length one. 
\end{rem}

\begin{exmp}
	\label{example2}
	Consider the graph $\mathcal{G} = (\mathcal{V},\mathcal{E})$ in Figure \ref{fig:graph2}. Moreover, consider the subsets of vertices $\mathcal{V}_1 := \{2,3\}$ and $\mathcal{V}_2 := \{6,7,8\}$.
		\begin{figure}[h!]
			\centering
			\scalebox{0.8}{
				\begin{tikzpicture}[scale=1]
				\node[draw, style=wball, label={180:$1$}] (1) at (0,0) {};
				\node[draw, style=wball, label={90:$2$}] (2) at (2,1) {};
				\node[draw, style=wball, label={-90:$3$}] (3) at (2,-1) {};
				\node[draw, style=wball, label={90:$4$}] (4) at (4,1) {};
				\node[draw, style=wball, label={-90:$5$}] (5) at (4,-1) {};
				\node[draw, style=wball, label={0:$6$}] (6) at (6,1.5) {};
				\node[draw, style=wball, label={0:$7$}] (7) at (6,0) {};
				\node[draw, style=wball, label={0:$8$}] (8) at (6,-1.5) {};
				\draw[-latex] (1) -- (2);
				\draw[-latex] (1) -- (3);
				\draw[-latex] (2) -- (4);
				\draw[-latex] (3) -- (4);
				\draw[-latex] (3) -- (5);
				\draw[-latex] (4) -- (6);
				\draw[-latex] (4) -- (7);
				\draw[-latex] (4) -- (8);
				\draw[-latex] (5) -- (7);
				\draw[-latex] (5) -- (8);
				\end{tikzpicture}
			}
			\caption{Graph used in Example \ref{example2}.}
			\label{fig:graph2}
		\end{figure}
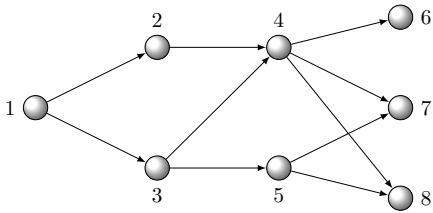 
		Clearly, the paths $\{(2,4), (4,6)\}$ and $\{(3,5), (5,7)\}$ form a set of two vertex-disjoint paths from $\mathcal{V}_1$ to $\mathcal{V}_2$. In fact, this set of vertex-disjoint paths is constrained, since there does not exist another set of two vertex-disjoint paths from $\bar{\mathcal{V}}_1 = \{2,3\}$ to $\bar{\mathcal{V}}_2 = \{6,7\}$. Note that there are also other sets of vertex-disjoint paths from $\mathcal{V}_1$ to $\mathcal{V}_2$. For example, also the paths $\{(2,4),(4,7)\}$ and $\{(3,5),(5,8)\}$ form a set of two vertex-disjoint paths. However, this set is \emph{not} constrained. To see this, note that there is another set of vertex-disjoint paths from $\bar{\mathcal{V}}_1 = \{2,3\}$ to $\bar{\mathcal{V}}_2 = \{7,8\}$, namely the set consisting of the paths $\{(2,4),(4,8)\}$ and $\{(3,5),(5,7)\}$. $\hfill{} \blacksquare$
\end{exmp}

The following theorem gives graph-theoretic conditions for identifiability of $(i,\mathcal{N}_i)$ in terms of constrained vertex-disjoint paths.

\begin{thm}
	\label{maintheorem}
	Consider a directed graph $\mathcal{G} = (\mathcal{V},\mathcal{E})$. Let $i \in \mathcal{V}$ and $\mathcal{C} \subseteq \mathcal{V}$. Then, $(i,\mathcal{N}_i)$ is identifiable from $\mathcal{C}$ if there exists a constrained set of $| \mathcal{N}_i |$ vertex-disjoint paths from $\mathcal{N}_i$ to $\mathcal{C}$.
\end{thm}

\begin{pf}
	Suppose that there exists a constrained set of $| \mathcal{N}_i |$ vertex-disjoint paths from $\mathcal{N}_i$ to $\mathcal{C}$. This means that there exists a subset $\bar{\mathcal{C}} \subseteq \mathcal{C}$ with $| \mathcal{N}_i | = |\bar{\mathcal{C}}|$ such that there is a constrained set of $| \mathcal{N}_i |$ vertex-disjoint paths from $\mathcal{N}_i$ to $\bar{\mathcal{C}}$. We want to prove that $(i,\mathcal{N}_i)$ is identifiable from $\bar{\mathcal{C}}$. By Lemmas \ref{lemma1} and \ref{lemma2}, this is equivalent to proving that the determinant of $(I-G(z))_{\mathcal{N}_i^c, \bar{\mathcal{C}}_{\phantom{i}}^c}$ is nonzero for all admissible network matrices $G(z)$. 
	
	We partition the vertex set $\mathcal{V}$ into four disjoint subsets, namely $\mathcal{R}, \mathcal{N}_i \setminus \bar{\mathcal{C}}$, $\mathcal{N}_i \cap \bar{\mathcal{C}}$, and $\bar{\mathcal{C}} \setminus \mathcal{N}_i$. Here the set $\mathcal{R}$ contains all the vertices that are not contained in any of the three other sets. Let $G(z)$ be an admissible network matrix. We compute
	\begin{equation*}
	M := (I-G)_{\mathcal{N}_i^c, \bar{\mathcal{C}}_{\phantom{i}}^c} = \begin{pmatrix}
	I-G_{\mathcal{R}, \mathcal{R}} & -G_{\mathcal{R}, \mathcal{N}_i\setminus\bar{\mathcal{C}}} \\
	-G_{\bar{\mathcal{C}}\setminus\mathcal{N}_i, \mathcal{R}} &  
	-G_{\bar{\mathcal{C}}\setminus\mathcal{N}_i, \mathcal{N}_i\setminus\bar{\mathcal{C}}} 
	\end{pmatrix},
	\end{equation*}
	where we have omitted the dependence on $z$ for the sake of brevity. Let $p := n - |\mathcal{N}_i|$ be the number of rows (and columns) of $M$. With $M$, we associate a weighted directed graph $\mathcal{G}_M = (\mathcal{V}_M,\mathcal{E}_M)$, where $\mathcal{V}_M = \{1,2,\dots,p\}$ and $\mathcal{E}_M := \{(k,l) \mid M_{lk} \neq 0 \}$. Furthermore, each edge $(k,l) \in \mathcal{E}_M$ is weighted by $M_{lk}$. Note that the graph $\mathcal{G}_M$ contains self-loops even though $\mathcal{G}$ was assumed to be simple. It is known that the determinant of $M$ can be expressed as a sum of the weights of spanning cycle families in $\mathcal{G}_M$. Recall from Section \ref{sectiongraphtheory} that a spanning cycle family in $\mathcal{G}_M$ is a collection of vertex-disjoint cycles such that each vertex in $\mathcal{V}_M$ appears in one of the cycles. 
	To be precise, we express $\det M$ as (cf. Theorem 3.1 of [\cite{Chen1971}])
	\begin{equation*}
	\det M = \pm \sum_\mathcal{F} (-1)^{N_\mathcal{F}} w(\mathcal{F}),
	\end{equation*}
	where $\mathcal{F}$ is a spanning cycle family in $\mathcal{G}_M$, $w(\mathcal{F})$ denotes the weight of the spanning cycle family (i.e., the product of the weights of all cycles in $\mathcal{F}$), the integer $N_\mathcal{F} \in \mathbb{N}$ denotes the number of cycles in $\mathcal{F}$, and the sum is taken over all spanning cycle families in $\mathcal{G}_M$. By our hypothesis, there exists a constrained set of $|\mathcal{N}_i|$ vertex-disjoint paths from $\mathcal{N}_i$ to $\bar{\mathcal{C}}$. By definition, this implies that there exists a constrained set of $|\mathcal{N}_i \setminus \bar{\mathcal{C}}|$ vertex-disjoint paths from $\mathcal{N}_i \setminus \bar{\mathcal{C}}$ to $\bar{\mathcal{C}} \setminus \mathcal{N}_i$. We denote this set of paths by $\mathcal{P}_1$. Let $\mathcal{N}_i \setminus \bar{\mathcal{C}} = \{v_1,v_2,\dots,v_r\}$ and $\bar{\mathcal{C}} \setminus \mathcal{N}_i = \{w_1,w_2,\dots,w_r\}$, where $r := |\mathcal{N}_i \setminus \bar{\mathcal{C}}|$. Without loss of generality, we can order the vertices in $\mathcal{N}_i \setminus \bar{\mathcal{C}}$ and $\bar{\mathcal{C}} \setminus \mathcal{N}_i$ such that there is a path in $\mathcal{P}_1$ from $v_j$ to $w_j$ for all $j = 1,2,\dots,r$. In terms of the graph $\mathcal{G}_M$, this means that there is a set of $r$ vertex-disjoint \emph{cycles} in $\mathcal{G}_M$ containing the nodes $\{|\mathcal{R}|+1,\dots,p\} \subseteq \mathcal{V}_M$. This is because \emph{columns} in $M$ corresponding to $\mathcal{N}_i \setminus \bar{\mathcal{C}}$ have the same indices as the \emph{rows} of $M$ corresponding to $\bar{\mathcal{C}} \setminus \mathcal{N}_i$. 
	We will denote this set of vertex-disjoint cycles by $\mathcal{F}_1$. The previous discussion implies that there exists a spanning cycle family in $\mathcal{G}_M$. Indeed, since the nodes in $\mathcal{V}_M$ that are not contained in $\mathcal{F}_1$ have a self-loop (with weight $1$), a spanning cycle family of $\mathcal{G}_M$ is given by the cycles in $\mathcal{F}_1$ together with self-loops on all remaining nodes in $\mathcal{V}_M$. We claim that \emph{every} spanning cycle family of $\mathcal{G}_M$ contains the cycles in $\mathcal{F}_1$. To see this, suppose on the contrary that there exists a spanning cycle family of $\mathcal{G}_M$ that does not contain all cycles of $\mathcal{F}_1$. Then, following the same reasoning as in [\cite{vanderWoude1991}] (see Section 3, page 37), we conclude that there exists a set of $r$ vertex-disjoint paths from $\mathcal{N}_i \setminus \bar{\mathcal{C}}$ to $\bar{\mathcal{C}} \setminus \mathcal{N}_i$ in $\mathcal{G}$ that is not equal to $\mathcal{P}_1$. However, this contradicts the hypothesis that $\mathcal{P}_1$ is \emph{constrained}. Therefore, each spanning cycle family in $\mathcal{G}_M$ contains the cycles in $\mathcal{F}_1$. Hence, we can rewrite the formula for $\det M$ as
	\begin{equation*}
	\det M = \pm w(\mathcal{F}_1) \sum_{\mathcal{F}_2} (-1)^{N_{\mathcal{F}_2}} w(\mathcal{F}_2).
	\end{equation*}
	Here $\mathcal{F}_2$ is a spanning cycle family of the subgraph of $\mathcal{G}_M$ obtained by removing all nodes (and incident edges) from $\mathcal{G}_M$ that appear in a cycle in $\mathcal{F}_1$. Moreover, $N_{\mathcal{F}_2}$ denotes the number of cycles in $\mathcal{F}_2$, and the sum is taken over all spanning cycle families of the subgraph of $\mathcal{G}_M$. Again, using Theorem 3.1 of [\cite{Chen1971}], we obtain
	\begin{equation*}
	\det M = \pm w(\mathcal{F}_1) \det (I-G_{\bar{\mathcal{R}}, \bar{\mathcal{R}}}),
	\end{equation*} 
	where $\bar{\mathcal{R}} \subseteq \mathcal{R}$ is the set of nodes in $\mathcal{R}$ that do not appear in one of the vertex-disjoint paths from $\mathcal{N}_i \setminus \bar{\mathcal{C}}$ to $\bar{\mathcal{C}} \setminus \mathcal{N}_i$ in $\mathcal{G}$. Finally, as $G$ is admissible, it satisfies property $P3$. Therefore $\det (I-G_{\bar{\mathcal{R}}, \bar{\mathcal{R}}}) \neq 0$. Moreover, since $w(\mathcal{F}_1)$ is the product of nonzero rational functions, $w(\mathcal{F}_1) \neq 0$. Therefore, $\det M \neq 0$. We conclude that $(i,\mathcal{N}_i)$ is identifiable from $\bar{\mathcal{C}}$ (and hence, from $\mathcal{C}$). $\hfill{} \square$
\end{pf}

\begin{exmp}
Consider the graph depicted in Figure \ref{fig:graph2}. For this example, let $\mathcal{C} := \{6,7,8\}$. Suppose we are interested in identifying the transfer functions associated with the edges from node $1$ to its out-neighbours $\mathcal{N}_1 = \{2,3\}$. To check that $(1,\mathcal{N}_1)$ is identifiable from $\mathcal{C}$, we use Theorem \ref{maintheorem}. In Example \ref{example2}, we already saw that there exists a constrained set of two vertex-disjoint paths from $\mathcal{N}_1$ to $\mathcal{C}$, namely the set consisting of the paths $\{(2,4),(4,6)\}$ and $\{(3,5),(5,7)\}$. Therefore, we conclude by Theorem \ref{maintheorem} that $(1,\mathcal{N}_1)$ is identifiable. That is, for any admissible network matrix $G(z)$ associated with $\mathcal{G}$ in Figure \ref{fig:graph2}, we can uniquely identify the transfer functions $G_{21}(z)$ and $G_{31}(z)$. $\hfill{} \blacksquare$
\end{exmp}
The following result gives graph-theoretic conditions under which \emph{all} transfer functions in $G$ are identifiable.
\begin{thm}
\label{theoremidentifiableG}
Consider a directed graph $\mathcal{G} = (\mathcal{V},\mathcal{E})$, and let $\mathcal{C} \subseteq \mathcal{V}$. Then, $\mathcal{G}$ is identifiable from $\mathcal{C}$ if for each $i \in \mathcal{V}$ there exists a constrained set of $| \mathcal{N}_i |$ vertex-disjoint paths from $\mathcal{N}_i$ to $\mathcal{C}$.	
\end{thm}
 Using Theorem \ref{theoremidentifiableG} we can show that all transfer functions appearing in the network of Figure \ref{fig:graph2} are identifiable, by measuring the node signals of just nodes $6,7,$ and $8$. 

\section{Conclusions}
\label{sectionconclusions}
In this paper we have considered the problem of identifiability of dynamical networks with partial node measurements. Unlike previous work [\cite{Hendrickx2018}] that considers \emph{generic} identifiability, we have considered identifiability \emph{for all} network matrices associated with the graph. We have introduced the new notion of \emph{constrained vertex-disjoint paths}. As our main result we have given a sufficient graph-theoretic condition for identifiability in terms of such paths. The authors are currently investiga\-ting necessary \emph{and} sufficient conditions for identifiability.

\bibliography{MyRef}             % bib file to produce the bibliography
                                                     % with bibtex (preferred)

\end{document}